\documentclass[12pt]{article}
\usepackage{amssymb}
\usepackage{epsfig}
\usepackage{epsf}
\usepackage{float}

\makeatletter
 
\@addtoreset{equation}{section}

\makeatother

\newtheorem{thm}[equation]{Theorem}
\newtheorem{lemma}[equation]{Lemma}

\newtheorem{propos}[equation]{Proposition}

\def\reff#1{(\ref{#1})}

\begin{document}

\def\E{{\mathbb E}}
\def\P{{\mathbb P}}
\def\R{{\mathbb R}}
\def\Z{{\mathbb Z}}
\def\V{{\mathbb V}}
\def\N{{\mathbb N}}
\def\X{{\cal X}}
\def\W{{\cal W}}
\def\G{{\cal G}}
\def\T{{\cal T}}
\def\I{{\cal I}}
\def\TT{\bar{{\cal T}}}
\def\II{\bar{{\cal I}}}
\def\C{{\C}}
\def\C{{\cal D}}
\def\n{{\bf n}}
\def\m{{\bf m}}
\def\b{{\bf b}}
\def\Var{{\hbox{Var}}}
\def\Cov{{\hbox{Cov}}}

\def\sqr{\vcenter{
         \hrule height.1mm
         \hbox{\vrule width.1mm height2.2mm\kern2.18mm\vrule width.1mm}
         \hrule height.1mm}}                  
\def\square{\ifmmode\sqr\else{$\sqr$}\fi}
\def\one{{\bf 1}\hskip-.5mm}
\def\limn{\lim_{N\to\infty}}
\def\given{\ \vert \ }
\def\ze{{\zeta}}
\def\be{{\beta}}
\def\la{{\lambda}}
\def\ga{{\gamma}}
\def\th{{\theta}}
\def\proof{\noindent{\bf Proof. }}
\def\A{{\bf A}}
\def\B{{\bf B}}
\def\C{{\bf C}}
\def\D{{\bf D}}
\def\MM{{\bf m}}
\def\w{\bar{w}}
\def\lnt{{\Lambda^N}}
\def\dlnt{\delta\Lambda^N_t}
\def\lno{\Lambda^N_0}
\def\dlno{\delta\Lambda^N_0}

\title{Asymptotic behavior of a stationary 
silo with absorbing walls}

\author{
{\large S. R. M. Barros, P. A. Ferrari} \\
{\large {\it Universidade de S\~{a}o Paulo}}\\
{\large N. L. Garcia},\\ {\it Universidade Estadual de Campinas}\\
{\large S. Mart\'{\i}nez},\\ {\it Universidad de Chile}
}
\date{}
\maketitle
\noindent {\bf Abstract}
We study the nearest neighbors one dimensional uniform q-model of force
fluctuations in bead packs [Coppersmith et al (1996)], a stochastic model to
simulate the stress of granular media in two dimensional silos. The vertical
coordinate plays the role of time, and the horizontal coordinate the role of
space. The process is a discrete time Markov process with state space
$\R^{\{1,\dots,N\}}$. At each layer (time), the weight supported by each grain
is a random variable of mean one (its own weight) plus the sum of random
fractions of the weights supported by the nearest neighboring grains at the
previous layer. The fraction of the weight given to the right neighbor of the
successive layer is a uniform random variable in $[0,1]$ independent of
everything. The remaining weight is given to the left neighbor. In the
boundaries, a uniform fraction of the weight leans on the wall of the silo.
This corresponds to \emph{absorbing boundary conditions}.  For this model we
show that there exists a unique invariant measure. The mean weight at site $i$
under the invariant measure is $i(N+1-i)$; we prove that its variance is
$\frac12(i(N+1-i))^2 + O(N^3)$ and the covariances between grains $i\neq j$
are of order $O(N^3)$. Moreover, as $N\to\infty$, the law under the invariant
measure of the weights divided by $N^2$ around site (integer part of) $rN$,
$r\in (0,1)$, converges to a product of gamma distributions with parameters
$2$ and $2(r(1-r))^{-1}$ (sum of two exponentials of mean $r(1-r)/2$). Liu
{\it et al} (1995) proved that for a silo with infinitely many weightless
grains, any product of gamma distributions with parameters $2$ and $2/\rho$
with $\rho\in [0,\infty)$ are invariant. Our result shows that as the silo
grows, the model selects exactly one of these Gamma's at each macroscopic
place.

\vskip 3mm
 {\bf Key words:} granular media, silos stress.

 {\bf AMS Classification:} Primary: 60K35, 82C20, 82C22

\vskip 3mm

\section{Introduction}

In the nearest neighbors one dimensional uniform q-model of force fluctuations
in bead packs the grains are arranged in layers inside a silo; each grain has
random weight with mean 1 and supports part of the weight of the grains of
higher layers.  The weight supported by each grain (plus its own weight) is
randomly distributed between the two neighboring grains of the following lower
layer.  The fraction of the weight given to the right nearest neighbor is
chosen uniformly in the interval $[0,1]$; the left nearest neighbor receives
the complementary weight. The walls of the silo absorb a random uniform
fraction of the weight supported by the boundary grains. The resulting model
is described by the following evolution equations. Fix a silo width $N$ and
for $i\in \{1,\dots,N\}$:
\begin{eqnarray}
W^N_t(i) &=& V_t(i) + W^N_{t-1}(i+1)U_{t-1}(i+1) 
+ W^N_{t-1}(i-1) (1 -
U_{t-1}(i-1))\,,\quad
\label{xx1}\\
W^N_t(0) &=& W^N_t(N+1) \;\equiv\; 0\,,
\label{xx2}
\end{eqnarray}
where $V_t(i)$ is the random weight of grain $i$ of level $t$ and $U_t(i)$ is
the random fraction of the weight grain $i$ of layer $t$ gives to grain $i-1$
of layer $t+1$. The boundary conditions \reff{xx2} represent the absorbing
wall. The random weights $V_t(i)$ have mean one: $\E V_t(i)=1$. The process
$W^N_t$ is Markov in the state space $(\R^+)^N$. We are interested in
stationary measures $\mu^N$ for this process. Under a stationary measure, the
mean weights $w^N(i)$ must satisfy the following equations. For $i\in
\{1,\dots,N\}$:
\begin{eqnarray}
w^N(i) &=& 1 + \frac12\, w^N(i-1) + \frac12\, w^N(i+1)\,,
\nonumber \\
w^N(0)&=&w^N(N+1)\;=\;0 \,.\label{1mu(i)}
\end{eqnarray}
obtained by taking expectations in \reff{xx1} and assuming that the law at
time $t$ is the same as the law at time $t-1$. These equations have a unique
solution, that turns out to be quadratic:
\begin{equation}
  \label{xx3}
  w^N(i)\;=\; i(N+1-i)\,.
\end{equation}
This implies that all candidates to be stationary measures must have the same
mean profile \reff{xx3}. Our first result uses \reff{xx3} to show the
existence of an invariant measure. Indeed, starting with any measure with this
quadratic profile, the means will be the same at all times. This guarantees
that the sequence of measures indexed by times $t\in\{0,1,2,\dots\}$ is tight
(\emph{i.e.}{} the mass does not escape to infinity; see the formal definition
after display \reff{tig} below) and by taking Ces\`aro limits along convergent
subsequences, one can show that there exists a stationary measure.  To prove
uniqueness we realize simultaneously two versions of the process starting with
different initial invariant distributions using the same sequence $U_t$ for
both evolutions. This is called \emph{coupling} in the probabilistic
literature. Under the coupling all weight added to the system after time zero
evolves identically in both versions, while the discrepant weight present at
time zero is distributed among the neighbors to eventually get lost at the
boundaries. As a consequence, both versions will have the same distribution in
the limit as $t\to\infty$. This is the key ingredient to show that there is a
\emph{unique} invariant measure $\mu^N$ for the silo.

Next step is to study properties of the invariant measure $\mu^N$. It seems
hard to explicitly describe the invariant measure for each finite $N$, but it
is possible to study the asymptotic behavior when $N$ grows. We show that the
covariances between the weights supported by different grains under the
invariant measure are bounded by a constant times $N^3$ and that the variance
at site $i$ differs from $(w^N(i))^2/2$ by at most a constant times $N^3$.
This implies that if one normalizes the stationary weights $W^N$ (distributed
according to $\mu^N$) dividing by $N^2$, one obtains, in the limit when
$N\to\infty$, that the covariances of the normalized weights vanish and that
the variance of $W^N([rN])/N^2$ converges to $(r(1-r))^2/2$. Here $[rN]$ is
the integer part of $rN$.

We also show that the law of the vector $W^N/N^2$ around site $[rN]$ converges
to a product of gamma distributions with parameters $2$ and $2/(r(1-r))$, the
sum of two independent exponentials of mean $r(1-r)/2$ each. To understand
this result divide \reff{xx1} by $N^2$ and take $N$ to infinity to get that
the limiting values must satisfy the same equation but in the infinite lattice
$\Z$ and without the $V_t$ (weightless grains). The limiting equations
correspond to the evolution of a model that we call \emph{infinite silo model}
(ISM). As Coppersmith et al (1996) did for the periodic case, a computation
shows that any space-homogeneous product of gamma distributions is invariant
for the ISM. The key observation to show this is the fact that if $X$ and $Y$
are independent exponentials with mean $\rho$ and $U$ is independent of $X,Y$
and uniformly distributed in $[0,1]$, then $U(X+Y)$ and $(1-U)(X+Y)$ are again
independent exponentials with mean $\rho$.  Actually we show more: following
an approach of Liggett (1976) we prove that \emph{all} invariant measures for
the ISM that are also shift invariant are mixtures of products of gammas with
parameters $2$ and $\rho$, with $\rho\in\R^+$. Hence, the distribution of the
vector $(W^N([rN]+\ell)/N^2:\ell\in\Z)$ converges to a mixture of products of
gamma distributions, as $N\to\infty$. Since the expected value of
$W^N([rN]+\ell)/N^2$ converges to $r(1-r)$ for all $\ell\in\Z$, it remains to
identify its limiting distribution as the gamma with parameters $2$ and
$2/(r(1-r))$. This is done using the fact that the normalized covariances
vanish as $N\to\infty$, as a consequence of our previous computation of the
limiting covariances.

Several authors proposed mathematical models for silos. The model in which a
grain lies its weight on to the lower neighbors was introduced by Harr (1977)
and explored by others, for example, Liu {\it et al.}  (1995) introduced the
model defined by \reff{W(i,t)}.  The model with zero boundary condition was
studied by Peralta-Fabi, M\'alaga and Rechtman (1997). Coppersmith {\it et
  al.}  (1996) developed mean field computations for these models and in the
uniform case, they conclude that the product of gamma distributions with
parameters $2$ and $2/\rho$ are invariant for a silo with periodic boundary
conditions and weightless grains; they also obtained analogous results for
dimensions $d\ge 2$.  The main contribution of Theorem \ref{111} is to prove
that as the size $N$ of the silo with zero boundary conditions increases, the
invariant measure around a macroscopic site $rN$ converges precisely to the
gamma distribution with parameters $2$ and $2/(r(1-r))$.  Rajesh and Majumdar
(2000) compute the space-space and space-time correlations for the periodic
model starting with the zero configuration. Claudin, Bouchaud, Cates and
Wittmer (1998) compute the spatiotemporal correlations in a continuum
approximation when the mass of each grain is rescaled so that it goes to zero
but Rajesh and Majumdar (2000) criticize the approach.  Socolar (1998) shows
numerical data for horizontal correlations for the model with periodic
boundary conditions.  Mueth, Jaeger and Nagel (1998) show experimental data
for force correlation in a silo filled with glass beads.  Essentially they
conclude, in agreement with the q-model, that the correlations vanish.
Lewandowska, Mathur and Yu (2000) computed the variances of stresses in the
$q$ model with weightless grains. Krug and Garcia (2000) and Rajesh and
Majumdar (2000) study related conserved-mass models with periodic boundary
conditions and show that in some cases the process has a stationary measure
with non zero correlations.

\section{Definitions and results}

To define the model let $(U_t):=\{U_t(i): i\in \Z,\, t\ge 0\}$ be a
family of independent uniform random variables in $[0,1]$ and
$(V_t):=\{V_t(i): i\in \Z,\, t\ge 0\}$ be a family of iid positive
random variables with mean 1 and variance $\alpha =\V V_t(i)<\infty$.
Furthermore assume $(V_t)$ and $(U_t)$ to be independent families.

Fix $N\ge1$, consider the finite box
\begin{eqnarray}
  \label{lnt}
  \Lambda^N&:=& \{1,...,N\}\nonumber
\end{eqnarray}
and denote $W^N_t(i)$ the weight carried by a grain located at
the $i$th position at level $t$. Fix an initial
configuration $W^N_0\in [0,\infty)^\lnt$ and define inductively 
\begin{eqnarray}
W^N_t(i) &=& V_t(i) + W^N_{t-1}(i+1)U_{t-1}(i+1) 
+ W^N_{t-1}(i-1) (1 -
U_{t-1}(i-1))\,,\quad
i\in \lnt\nonumber\\
W^N_t(0) &=& W^N_t(N+1) \;\equiv\; 0\,.
\label{W(i,t)}
\end{eqnarray}
Let $W^N_t= (W^N_t(i):i\in \lnt)$; then $(W^N_t: t\ge 1$) is a discrete time
Markov chain on $[0,\infty)^{\lnt}$.  Each grain $j$ of layer $t$ gives a
fraction chosen uniformly in $[0,1]$ of its own weight plus the total weight
it supports from the previous layers to grain $j-1$ of the successive layer
$t+1$ (which we can think is below $t$) and the remaining to grain $j+1$. The
weight distributed to grains outside $\lnt$ is thought of as being absorbed by
the walls of the silo at sites $0$ and $N+1$.

For a measure $\nu$ on $[0,\infty)^\lnt$ let $\nu S^N(t)$ be the measure
defined by
$$
\nu S^N(t) f \;=\;\E_\nu f(W^N_t) \;:= \;\int \nu(dW) \,\E(f(W^N_t)\,|\,W^N_0=W).
$$
where $\E$ and
$\P$ are the expectation and probability defined with respect to the
probability space induced by $(U_t:t\ge 0)$ and $(V_t:t\ge 0)$.

We say that a measure $\mu^N$ is \emph{invariant} for the process
$W^N_t$ if $\mu^NS^N(t)=\mu^N$.

If $W^N_t$ has an invariant measure $\mu^N$ its mean heights $w^N(i):=\int
\mu^N(dW)W(i)$ have to satisfy the following system of equations (taking
expectations in \reff{W(i,t)} and assuming that the distribution of $W^N_t$
does not depend on $t$):
\begin{eqnarray}
w^N(i) &=& 1 + \frac12 w^N(i-1) + \frac12 w^N(i+1)\,,
\mbox{ \ for $i\in \lnt$\,;}
\nonumber \\
w^N(i)&=&0,\;\;\;\;\mbox{ for }i\in \Z\setminus\lnt \,.\label{mu(i)}
\end{eqnarray}
Notice that the only finite solution $w^N(\cdot)$ of \reff{mu(i)} has quadratic
profile:
\begin{eqnarray}
  w^N(i) &=& i(N+1-i)\,, \mbox{ \ for $ i=0, \ldots, N+1$}\,.
\label{rw(i)} 
\end{eqnarray}
(It is the expected time to exit $\lnt$ for a symmetric nearest neighbors
random walk starting at $i$.) 

Our first result says that for each $N$ the silo
with absorbing boundary conditions admits a unique invariant measure and that
the process starting with any initial condition converges to the invariant
measure.

\begin{thm} 
\label{2.4}
There exists a unique invariant measure $\mu^N$ for the process $W^N_{t}$
defined by \reff{W(i,t)}. The mean weights under $\mu^N$ are finite and
satisfy \reff{rw(i)}. Furthermore, for any initial configuration $W^N_0$,
the law of $W^N_t$ converges to $\mu^N$.
\end{thm}

Denote the covariances of weights $i$ and $j$ under $\mu^N$ by
\begin{eqnarray}
\sigma^N(i,j) &:=& \int\mu^N(dW)\, W(i)W(j) \,-\, w^N(i)\,w^N(j)\,. \nonumber
\end{eqnarray}
We show the following bounds for the covariances.
\begin{thm}
  \label{z6}
  The covariances
  of $\mu^N$ are of order $N^3$: there exists a positive constant $C$ such
  that for all $N$
\begin{equation}
  \label{z7}
 |\sigma^N(i,j)|\;\le\; CN^3,\qquad i\neq j,\; i,j\in
\Lambda_N
\end{equation}
and 
\begin{equation}
  \label{z30}
  \Bigl|\sigma^N(i,i) - 
{(w^N(i))^2\over 2}\Bigr|\le CN^3 ,\qquad i\in \Lambda_N\,.
\end{equation}
\end{thm}
Numerical analysis performed in Section 2 suggests 
\begin{equation}
  \label{p1}
  \lim_{N\to\infty}{\sigma^N([xN],[yN])\over N^3}\; =\; -\,{r(x,y)\over
  3}\;;\qquad 
  \lim_{N\to\infty}{\sigma^N([xN],[xN]) -  {(w^N([xN]))^2\over 2}\over N^3} 
\;=\; -\,{r(x,x)\over 2}\,,
\end{equation}
for $x,y\in[0,1]$, $x\neq y$, where $r(x,y)$ is the unique weak solution of 
the PDE
\begin{eqnarray}
  \label{pde}
-\Delta u - a(x,y) \Bigl(\Delta u + 2\frac{\partial^2u}{\partial x\partial
  y}\Bigr)& =& 
f(x,y)\ \ \mbox{in} \ \  
\Omega=(0,1)\times(0,1) \nonumber\\
u & = & 0 \ \ \mbox{in} \ \ \partial\Omega\,,
\end{eqnarray}
where $a(x,x)=1/2$ and $a(x,y)=0$ for $x\neq y$ and $\Delta$ stands for the
Laplacian operator. The right-hand-side is given by $f(x,y) =3\sqrt{2} (6(x(1-x)-1)
\delta(x-y)$, where $\delta$ is Dirac's distribution.  An approximate plot
of the function $r$ obtained numerically can be seen in Figure 1. The
numerical solution also suggests that
\begin{eqnarray}
  \lim_{N\to\infty} {\sigma^N(i,j)\over N^2} &=& c(i,j)\; >\; 0\,. \label{p2}  
\end{eqnarray}
The meaning of \reff{p2} for $i=j=1$ is that the variance of the weight
discharged to the left wall at each level is of order $N^2$, the square of the
mean. Hence the fluctuations are of order $N$, the order of the mean. By
symmetry the same happens to the right wall. On the other hand, \reff{p2} for
$i\neq j$ means that the stationary correlations divided by the square of the
mean stay strictly positive around the origin; this discards the hypothesis
that $W^N(\cdot)/N^2$ converges to a product measure. In contrast, the
correlations divided by the square of the mean vanish in the interior of the
silo as we show in the next theorem.

For $\rho>0$
define
\begin{equation}
\label{nur}
\nu_{\rho}:= \hbox{\ \  product measure on }
[0,\infty)^{\Z} \hbox{\ \ with marginals Gamma }(2, 2/ \rho)
\end{equation}
 (the sum of two
independent exponentials with mean $\rho/2$ each). 
Our main result says that if $W$ is distributed according to the invariant
measure $\mu^N$, then as $N\to\infty$, the law of $W/N^2$ around the site
$[rN]$ converges to $\nu_{r(1-r)}$. Let $\tau_{i}$ be the translation operator
by $i$ on $\Z$ and $\Theta_k$ be the operator ``divide by $k$'': for $W\in
[0,\infty)^\Z$, $(\tau_{i}W)(j)= W(j-i)$ and $(\Theta_kW)(j) = W(j)/k$.

\begin{thm} \label{111} Let $\mu^N$ be the unique invariant measure for the
  process $W^N_{t}$. Then for each $r \in (0,1)$
\begin{equation}
\label{110}
\lim_{N\to\infty}\tau_{[rN]}\Theta_{N^2}\mu^N 
= \nu_{r(1-r)}
\end{equation}
weakly, where $[\cdot]$ is the integer part.
\end{thm}

The proof of Theorem \ref{111} lies on two results. The first one is related
with the limit as $N\to\infty$ of the weights of the silo model divided by
$N^2$ in sites around a macroscopic position $rN$. Its distribution converges
to the law of another process in $(\R^+)^\Z$ whose evolution is the same
described by \reff{W(i,t)} but with $V_t(i)\equiv 0$ ---no extra mass is added
in this system. This process is called the \emph{infinite silo model, ISM}. We
follow an approach of Liggett (1976) to show that in the nearest neighbors
uniform case, the invariant and translation invariant measures for the ISM are
mixtures of product of Gamma distributions (Proposition \ref{lemma3} below).
Combining this with the fact that the stationary space-space correlations go
to zero, proved in Theorem \ref{z6}, we can single out the product measure as
the limit in \reff{110}.

Theorems \ref{2.4} and \ref{z6} are proven in Section \ref{s2}. A numerical
study of the corrections for the correlations in Theorem \ref{z6} is given at
the end of Section \ref{s2}. Theorem \ref{111} is proven in Section 4. In
Section 3 we study the infinite silo model.

\section{The finite stationary silo}
\label{s2}

In this section we study the silo for fixed width $N$. We prove the existence
of an invariant measure and study its correlations. 

Notice that \reff{W(i,t)} describes the evolution of two independent systems:
$\{W^N_t(i)\,:\,i+t$ even$\}$ and $\{W^N_t(i)\,:\,i+t$ odd$\}$. It is
convenient to keep simultaneously track of both systems.

\paragraph{Proof of Theorem \ref{2.4}}
 If we consider any initial distribution $\nu_{0}$ for $W^N_0$,
not necessarily invariant, satisfying $\E_{\nu_0}[W^N_0(i)] = i(N+1-i)$,
$x\in \lno$, then for all $t\ge 0$
\begin{equation}
\label{1.1}
\E_{\nu_0}[W^N_t(i)] = i(N+1-i),\;\; i\in \lnt\,.
\end{equation}
Since by Markov inequality 
\begin{equation}
  \label{tig}
 \P_{\nu_0}[W^N_t(i) > M] \le \frac{\E_{\nu_0}[W^N_t(i)]}{M}\,,
\end{equation}
the sequence $( W_t: t\ge 0)$ is tight [this means that for all
$\varepsilon>0$ there exists a $K>0$ such that $\P(W_t(i)>K)<1-\varepsilon$
for all $i,t$]. Tightness implies that the Cesaro means $\nu_T:=(1/T)\sum_{t=1}^T
S^N(t)\nu_0$ converge through subsequences to an invariant measure (see
Liggett (1976), for instance) called~$\mu^N$. 

If $(t_k)$ is a subsequence such that $\nu_{t_k} \to \mu^N$, as $k \to\infty$,
then there exist $\overline W_k$ with law $\nu_{t_k}$ and $W$ with law
$\mu^N$, all defined in the same probability space, such that for
$i=1,\dots,N$, $\overline W_k(i) \to W(i)$ almost surely (Theorem 8.1, page 25
in Thorisson (2000); this is usually called ``almost sure Skorohod
representation for distribution convergence'').  Since the means are
independent of $t$, Fatou's Lemma ($\E (\liminf_k \overline W_k(i)) \le
\liminf_k \E \overline W_k(i)$) implies that $\E W(i) \le i(N+1-i)$, that is,
the means under $\mu^N$ are bounded by the quadratic profile~\reff{rw(i)}.
But if the means are bounded and satisfy the equations~\reff{mu(i)}, then
$\mu^N$ must have the quadratic profile~\reff{rw(i)}.

In order to prove uniqueness we use coupling. Construct two processes $X_t$
and $Y_t$ satisfying (\ref{W(i,t)}) using the same uniform random variables
$(U_t(i):i,t)$ and $(V_t(i):i,t)$.  Assume the initial configurations satisfy
\begin{equation}
  \label{x2}
X_0(i) \ge Y_0(i) \mbox{  for }i\in \lnt\,.
\end{equation}
Then it is easy to see that at any time $t\ge 0$ 
\begin{equation}
  \label{x3}
X_t(i) \ge Y_t(i) \mbox{  for }i\in \lnt
\end{equation}
(in this case we say that the process is \emph{attractive}).
Denote
\[ D(t) := \frac{1}{N} \sum_{i=1}^{N} (X_t(i) - Y_t(i))\,.\]
We have
\begin{equation}
D(t) = D(t-1) - \frac{1-U_{t-1}(1)}{N} [X_{t-1}(1) - Y_{t-1}(1)] 
- \frac{U_{t-1}(N)}{N}[X_{t-1}(N) -
Y_{t-1}(N) ]\,. \label{D(k)}
\end{equation}
The numbers $D(1), D(2), \ldots$ form a bounded below by $0$ non-increasing
sequence of nonnegative random variables. Hence $D(t)$ converges almost
surely. This implies that its increments must go almost surely to zero. This
and \reff{D(k)} imply that
\[
[X_{t-1}(1) - Y_{t-1}(1)] \rightarrow 0
\;\;\mbox{ and }\;\;
[X_{t-1}(N) - Y_{t-1}(N)] \rightarrow 0\,,
\]
almost surely, as $t \rightarrow \infty$. Looking at the equations
\[
X_{t-1}(1) - Y_{t-1}(1)  =  (1-U_{t-2}(2))\,[X_{t-2}(2) - Y_{t-2}(2)] \]
and 
\[ X_{t-1}(N) - Y_{t-1}(N) = U_{t-2}(N-1)\, [X_{t-2}(N-1) -
Y_{t-2}(N-1)]\,,\]
it is immediate to see that
\[
X_{t-2}(2) - Y_{t-2}(2) \rightarrow 0
\;\;\;\mbox{ and }\;\;\;
X_{t-2}(N-1) - Y_{t-2}(N-1)\,, \rightarrow 0
\]
almost surely, as $t \rightarrow \infty$. Using the same argument we
have that for all $\ i\in \lnt$
\[
X_t(i) - Y_t(i) \rightarrow 0, \;\;\;\mbox{ as }t \rightarrow \infty\,,
\]
almost surely as $t \rightarrow \infty$.

Let $(Z_0(i): i\in \lnt)$ and $(W_0(i): i\in \lnt)$ arbitrary random
vectors (not necessarily ordered). Define
\[ 
X_0(i) = Z_0(i)\vee  W_0(i)
\mbox{ and }
Y_0(i) = Z_0(i)\wedge  W_0(i)\,.\]
We have
\[ 
Y_0(i)\le W_0(i) \le X_0(i)\,.
\] 
Couple the four processes using equation (\ref{W(i,t)}). This means that the
four evolutions are performed using the same random variables $V_t(i)$
and $U_t(i)$. The previous argument implies
\[ 
W_t(i) - Y_t(i) \rightarrow 0
\mbox{ and }
X_t(i) - W_t(i) \rightarrow 0\,,
\]
almost surely as $t \rightarrow \infty$. Since
\[ 
W_t(i) - Y_t(i) \le |W_t(i) - Z_t(i) | \le
X_t(i) - W_t(i),
\]
we have
\begin{equation}
  \label{x4}
 W_t(i) - Z_t(i) \rightarrow 0\,, 
\end{equation}
almost surely as $t \rightarrow \infty$. 

Hence, if we pick $W_0$ from the distribution $\mu^N$ and $Z_0$ from an
arbitrary invariant distribution $\nu$, we have $\mu^N=\nu$. This proves that
$\mu^N$ is the unique invariant measure. The limit \reff{x4} also shows that
the process $W^N_t$ starting from an arbitrary configuration converges to
$\mu^N$. \square

\paragraph{Proof of Theorem \ref{z6}} 
At the end of the proof we show that $|\sigma^N(i,j)|<\infty$ for all $i,j$
and $N$.
If $W^N_0$ is distributed according to the invariant measure, so is
$W^N_1$, and one can use \reff{W(i,t)} to show that $\sigma^N(i,j)$ satisfies
the
system of equations:
\begin{eqnarray}
\sigma^N(i,i) & = & \alpha\,+\,
\frac{1}{3}\sigma^N(i+1,i+1)\,+\,\frac{1}{3}\sigma^N(i-1,i-1) \,+\,
\frac{1}{4}\sigma^N(i-1,i+1)\nonumber \\ 
&& \quad +\,\frac{1}{4}\sigma^N(i+1,i-1) \,+\,
\frac{1}{12}(w^N(i+1))^2\,+\,\frac{1}{12}(w^N(i-1))^2,\quad i\in\Lambda^N
\nonumber \\ 
\sigma^N(i,i+2) & = &
\frac{1}{4}\sigma^N(i+1,i+3)\,+\,\frac{1}{4}\sigma^N(i-1,i+3) \,+\,
\frac{1}{4}\sigma^N(i-1,i+1) \nonumber \\
&& \quad +\,\frac{1}{6}\sigma^N(i+1,i+1) \,-\,
\frac{1}{12}(w^N(i+1))^2  , \qquad i\in\{1,\dots,N-2\}\nonumber  \\ 
\sigma^N(i,i-2) & = &
\frac{1}{4}\sigma^N(i+1,i-1)\,+\,\frac{1}{4}\sigma^N(i-1,i-3) \,+\,
\frac{1}{4}\sigma^N(i+1,i-3) \nonumber \\
&& \quad +\,\frac{1}{6}\sigma^N(i-1,i-1) \,-\,
\frac{1}{12}(w^N(i-1))^2  , \qquad i\in\{3,\dots,N\}\nonumber  \\ 
\sigma^N(i,j) & = &
\frac{1}{4}\sigma^N(i+1,j+1)\,+\,\frac{1}{4}\sigma^N(i-1,j-1) \,+\,
\frac{1}{4}\sigma^N(i-1,j+1) \nonumber \\
&& \quad \,+\,\frac{1}{4}\sigma^N(i+1,j-1), 
\qquad i,j\in\Lambda^N,\, \quad |i- j| 
\ge 2 \nonumber \\  
\sigma^N(i,j) & = &0\,,\qquad i\in\{0,N+1\}\hbox { or }
j\in\{0,N+1\}\,.\label{cov} 
\end{eqnarray}
Define, for $i \in\{ 0,1,\ldots, N+1\}$ :
\begin{eqnarray}
R(i,i) & = & -2 \sigma^N(i,i) +(w^N(i))^2\,,
\nonumber \\
R(i,j) & = & - 3 \sigma^N(i,j) \quad\mbox{ for } i\neq j\,.\label{r}
\end{eqnarray}
The system \reff{cov} translates into the following system for
the matrix $R$:
\begin{eqnarray}
R{(i,j)} &=& K(i)\,\one\{i=j\} + \sum_{(i',j')} p((i,j),(i',j'))
R{(i',j')}\qquad i,j\in\Lambda^N\,, \nonumber\\
R{(i,j)} &=& 0 \qquad  i\in\{0,N+1\}\hbox { or }
j\in\{0,N+1\}\,,\label{tij}
\end{eqnarray}
where 
\begin{eqnarray}
  \label{s14}
  K(i)&=&(w^N(i))^2 - \frac{(w^N(i+1))^2}{2} -
  \frac{(w^N(i-1))^2}{2}-2\alpha\nonumber\\ 
&=&6i(N+1-i)-(N+1)^2-(1+2\alpha) \,,
\end{eqnarray}
for $i\in\Lambda^N$ and 
\begin{equation}
\label{z8}
p((i,j),(i',j')) := \left\{ \begin{array}{ll}
                          {1}/{4}, & \mbox{ if } i \neq j \mbox{
                          and } |i-i'|.|j-j'|=1 \\
                          {1}/{3}, & \mbox{ if } i = j, i'=j' \mbox{
and
                           } |i-i'|=1 \\
                          {1}/{6}, &\mbox{ if } i = j, i' \neq j'
                          \mbox{ and } |i-i'|.|j-j'|=1.\\
                          0 &\mbox{ otherwise.}
                          \end{array}
                   \right.
\end{equation}
for $i,j\in\{1,\dots,N\}$ and 
$$
p((i,j),(i',j')):=\cases {1 &if $(i',j')=(i,j)$\cr
0& otherwise}
$$
if $i\in\{0,N+1\}$ or $j\in\{0,N+1\}$ (this means absorbing borders). Let
$X^{(i,j)}_n$ be a Markov chain with transition probabilities $p(\cdot,\cdot)$
and initial state $(i,j)$. Defining
\begin{equation}
  \label{z10}
T^{(i,j)}:= \sum_{n\ge 0} \sum_{\ell=1}^N
\one\{X^{(i,j)}_n =(\ell,\ell)\}\, K(\ell) 
\end{equation}
we have that $\E T^{(i,j)}$ is the unique finite solution of \reff{tij}.
Assuming that $|\sigma^N(i,j)|$ is finite, then $|R(i,j)|$ is also finite and
$R(i,j)=\E T^{(i,j)}$.

To get an upperbound for $R{(i,j)}$ we consider a random walk $\tilde
X^{(i,j)}_n$ on the diamond 
$$
{\mathbf \Delta}^N:=\{(i,j)\in\Z^2: i+j \mbox{ even},\,2\le i+j\le 2N,\, | i-j|\le
N+1\} 
$$
with transition probabilities 
\begin{eqnarray}
  \label{x5}
 && \tilde p((i,j),(i',j')) \;=\; p((i,j),(i',j')) \mbox{ given by \reff{z8}}
 \nonumber\\
&&\qquad \mbox{ if } (i,j)\in
  \{(i,j)\in\Lambda^N\times\Lambda^N: 4\le i+j\le 2N-4,\,  
|i-j|\le N-1\}
\end{eqnarray}
(the interior of the diamond); periodic boundary conditions along two of the
sides of the diamond:
\begin{eqnarray}
  \label{z9}
  \tilde p((i,j),(N+1-j,N+1-i)) &=& \tilde p((N+1-j,N+1-i),(i,j)) \;=\;1/4
  \nonumber\\ 
\lefteqn{ \hbox{ if
    }\quad i+j = 2,\quad i\neq j} \nonumber
\end{eqnarray} 
\begin{eqnarray}
  \label{z9a}\tilde p((i,j),(i',j')) &=& \tilde
  p((N+1-j,N+1-i),(N+1-j',N+1-i')) \;=\;1/4 
 \nonumber\\  
\lefteqn{\hbox{ if
    }\quad i+j = 2,\quad i\neq j,\quad |i-i'|.|i-j'|=1,\quad (i',j')\in
 {\mathbf \Delta}^N }\nonumber
\end{eqnarray} 
\begin{eqnarray}
\tilde p((1,1),(i,j)) &=& \cases{
1/3&if $(i,j)\in\{(2,2),(N,N)\}$\cr
1/6&if $(i,j)\in\{ (2,0),(0,2)\}$}\\
\tilde p((N,N),(i,j)) &=& \cases{
1/3&if $(i,j)\in\{ (1,1),(N-1,N-1)\}$\cr
1/6&if $(i,j)\in\{ (N-1,N+1),(N+1,N-1)\}$}
\end{eqnarray}
and absorbing probabilities along the other two sides:
\begin{equation}
  \label{z33}
  \tilde p((i,j),(i,j))\;=\;1\,,\qquad (i,j)\in{\mathbf \Delta}^N,\; |i-j|= N+1\,.
\end{equation}
The point of this random walk is that it is ``one dimensional'' in the sense
that the one-dimensional process defined by $Z^{i-j}_n:=(\tilde
X^{(i,j)}_n)_1-(\tilde X^{(i,j)}_n)_2$ (the difference of the coordinates of
$\tilde X^{(i,j)}_n$) is Markovian and has transition probabilities
\begin{equation}
  \label{z34}
  p_0(i,j) = \cases {1/4&if  $i\notin\{-N-1,0,N+1\}$ and $|i-j|=2$\cr
1/2 &if  $j=i$ and $i\notin\{-N-1,0,N+1\}$ \cr
2/3 &if $i=j=0$\cr 
1/6 &if $i=0$ and and $|i-j|=2$\cr 
1&if $i=j$ and $i\in\{-N-1,N+1\}$}
\end{equation}
for even $i\in \{-N-1,\dots,N+1\}$.

Define $ \overline K= \max_\ell K(\ell)$, $\underline K = \min_\ell K(\ell)$. 
We have $\underline K< 0< \overline K$ and both
$\underline K, \overline K$ are of order of $N^2$. Define
\begin{eqnarray}
  \label{s1}
 \overline T^{(i,j)}&:=& \sum_{n\ge 0} \sum_{\ell=1}^N
\one\{\tilde X^{(i,j)}_n =(\ell,\ell)\} \overline K\nonumber\\
  \underline T^{(i,j)}&:=& \sum_{n\ge 0} \sum_{\ell=1}^N
\one\{\tilde X^{(i,j)}_n =(\ell,\ell)\} \underline K\,.\nonumber
\end{eqnarray}
A simple coupling argument shows that
for all $i,j\in\Lambda^N$,
\begin{equation}
  \label{z11}
  \underline T^{(i,i)}\; \le\;\underline T^{(i,j)}\; \le\; T^{(i,j)}\;
  \le\;\overline T^{(i,j)}\; \le\;\overline T^{(i,i)}\,. 
\end{equation}
On the other hand $\underline T^{(i,i)}$ and $\overline T^{(i,i)}$ do not
depend on $i$ and can be expressed in function of a geometric random sum of
geometric random variables:
\begin{equation}
  \label{z12}
  \underline T^{(i,i)} = \underline K \sum_{k=1}^S M_k\;;\qquad
\overline T^{(i,i)} = \overline K \sum_{k=1}^S M_k\,,
\end{equation}
where $\P(S\ge n) = ({N-1\over N+1})^{n-1}$, $n\ge 1$ and $\P(M_k\ge m) =
(2/3)^{m-1}$, $m\ge 1$. Furthermore $(M_k)$ are iid and independent of $S$.
The random variable $S$ counts the number of times the chain $Z_n$ starting at
$2$ or $-2$ comes back to zero before arriving to $N+1$ or $-N-1$ while $M_k$
is the time the chain $Z_n$ spends in its $k$th visit to the origin.  Taking
expectations in \reff{z12} and using Wald identity we get
\begin{equation}
  \label{z13}
  \E \overline T^{(i,i)} = {3\over 2} (N+1)\,\overline K\;;
\qquad \E \underline T^{(i,i)} = {3\over 2} (N+1)\,\underline K\,.
\end{equation}
Taking expectations in \reff{z11} and using the fact that both $\underline K$
and $\overline K$ are of order $N^2$, we have proved that if $R(i,j)$ is
finite,
then
\begin{equation}
  \label{z333}
  -C N^3\le \underline R(i,j)\le R(i,j)\le \overline R(i,j)\le CN^3\,.
\end{equation}
This, \reff{r} and \reff{1.1} show \reff{z7} and \reff{z30}.

\paragraph{Finiteness of $\sigma^N(i,j)$}
Let $W^N_0$ be a non negative random configuration chosen according to a
product measure with marginals satisfying
\begin{eqnarray}
  \label{z18}
\E W^N_0(i) &=& w^N(i)\nonumber\\
  \V W^N_0(i) &=& (w^N(i))^2/2
\end{eqnarray}
(the covariances are null). The covariances at time $t$ are given by
\begin{equation}
  \label{z19}
   \sigma^N_t(i,j) = \cases{\frac12(w^N(i))^2 -
   \frac12R^0_t(i,i) & if   $i=j$\cr
-  \frac13R^0_t(i,j)& if
  $i\neq j$}\,,
\end{equation}
where $ R^0_t(i,j)$ evolves according to 
\begin{equation}
 R^0_t{(i,j)} = \one\{i=j\}K(i) + \sum_{(i',j')} p((i,j),(i',j'))
 R^0_{t-1}{(i',j')}\,, \label{tij1}
\end{equation}
with $ R^0_0{(i,j)}=0$ for $i,j\in\{0,\dots,N+1\}$. Let $\overline R_t$ and
$\underline R^0_t$ be the systems defined in the diamond ${\mathbf \Delta}^N$
with initial condition $\overline R^0_0(i,j)=\underline R^0_0(i,j)=0$,
$i,j\in{\mathbf \Delta}^N$ and evolving with the equations
\begin{eqnarray}
  \label{s3}
\overline R^0_t{(i,j)} &=& \one\{i=j\}\overline K + \sum_{(i',j')}
 \tilde p((i,j),(i',j')) 
 \overline R^0_{t-1}{(i',j')}. \label{s4}\\
\underline R^0_t{(i,j)} &=& \one\{i=j\}\underline K + \sum_{(i',j')}
 \tilde p((i,j),(i',j')) 
\underline R^0_{t-1}{(i',j')}\,. \label{s5}  
\end{eqnarray}
It is clear that 
\begin{eqnarray}
  \label{z20}
   \underline R^0_t(i,j) &\le&  R^0_t (i,j)
\;\le\; \overline R^0_t(i,j)\,.\nonumber
\end{eqnarray}
On the other hand
\begin{eqnarray}
  \label{zz20}
 \overline R^0_t(i,j)\;\le\; \overline
R_{t}(i,j) \;\le \; CN^3\,,
\end{eqnarray}
where $\overline R_{t}\equiv\overline R$ is the stationary solution of
\reff{s4}.  That is, $\overline R^0_t(i,j)$ is bounded by the stationary
solution of \reff{s4} which in turn is bounded as in \reff{z333}.  The same
argument shows that $|\underline R^0_t(i,j)|$ is uniformly bounded in $t$.
Hence $|R^0_t{(i,j)}|$ is uniformly bounded in $t$ and by \reff{z19} so is
$|\sigma^N_t(i,j)|$. To show the finiteness of $\sigma^N(i,j)$; we need to
take the limit as $t\to\infty$. Since the distribution $\nu^N_t$ of
$W_t^N(i)$, converges weakly to $\mu^N$, as $t \to \infty$, we can again use
Theorem 8.1, page 25 in Thorisson (2000) to define $W_t^N(i)$ and $W^N(i)$ in
the same space with marginal laws $\nu^N_t$ and $\mu^N$, respectively, such
that $W_t^N(i)$ converges almost surely to $W^N(i)$, as $t\to\infty$. Now,
$\sigma^N(i,i) = \V W^N(i) = \E (W^N(i))^2 - (\E W^N(i))^2$. Since $\E W^N(i)=
\E W_t^N(i)$ for all $t$, Fatou's Lemma (which in this case says $\E
(\liminf_t (W_t^N(i))^2) \le \liminf_t \E (W^N(i))^2$) implies
\begin{eqnarray}
  \label{z21}
 \sigma^N(i,i) &\le& \liminf_{t\to \infty}
 \sigma^N_t(i,i)\;\le\; 
 \displaystyle{(w^N(i))^2\over 2}+CN^3 
\end{eqnarray}
by \reff{z19} and \reff{z333}. This settles the finiteness of the diagonal
covariances. By \reff{cov}, also the other covariances must be finite.
\square

\begin{figure}[ht]
\begin{center}
\epsfig{figure=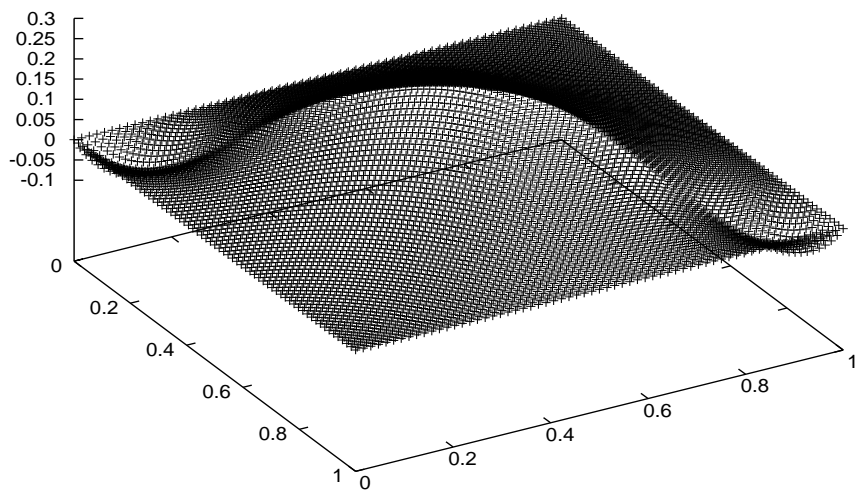,height=14cm,width=9cm,angle=-90}
\caption{ Graphic of the computed solution $r=R/(N+1)^3$ for
$N=127$.}
\end{center}
\end{figure}

\begin{figure}[ht]
\begin{center}
\epsfig{figure=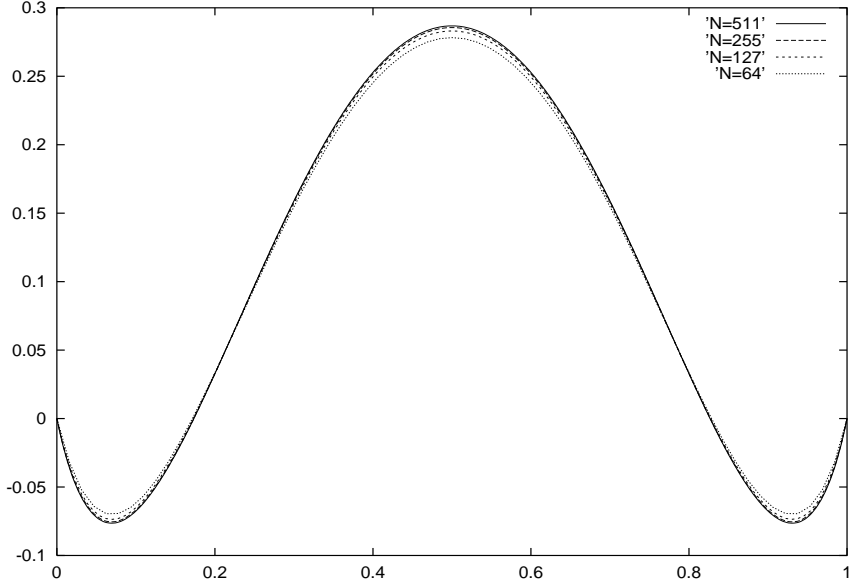,height=12cm,width=8cm,angle=-90}
\caption{ Graphic of the computed solution $r=R/(N+1)^3$ over the
  diagonal ($i=j$) for several values of $N$.}
\end{center}
\end{figure}

\paragraph{Numerical Analysis of the equation for the correction $R$}
We carried out a numerical investigation about the solution of equation
\reff{tij}. This equation can be seen as a discretized
version of the partial differential equation \reff{pde} in the following way.
We first divide system \reff{tij} by $(N+1)^3$, obtaining a system for the new
variable $r_{i,j}:=R_{i,j}/(N+1)^3$, which is to be seen as an approximation
to the value of a function $r(x,y)$ defined on $\Omega=[0,1]\times [0,1]$ at
the point $(x_i,y_j)=(ih,jh)$, where $h=1/(N+1)^3$. The equations for $i\neq
j$ are of the form
$$
4r_{i,j}-r_{i-1,j-1}-r_{i-1,j+1}-r_{i+1,j-1}-r_{i+1,j+1}=0
$$
which is a discrete form of $\Delta r=0$, where $\Delta$ denotes the
(rotation invariant) Laplace operator. Now, for $i=j$ we multiply \reff{tij}
by $3h=3/(h^2(N+1)^3)$ and obtain the equations:
\begin{eqnarray}
\frac{1}{2h^2}\left(
    4r_{i,i}-r_{i-1,i-1}-r_{i-1,i+1}-r_{i+1,i-1}-r_{i+1,i+1} 
\right)+ \frac{1}{2h^2}\left( 2r_{i,i}-r_{i-1,i-1}-r_{i+1,i+1} \right) 
\nonumber \\
=\; 3\sqrt{2}\left( 6x_i(1-x_i)-1-(1+2\alpha)h^2 \right)
\delta_{\sqrt{2}h}(i-i)\,,\qquad \label{discpde}
\end{eqnarray}
where $x_i=ih$ and $\delta_{\sqrt{2}h}$ is a piecewise linear approximation to
Dirac's distribution, assuming the value $1/(\sqrt{2}h)$ at the origin, and
being zero at the other points of an equally spaced mesh (with mesh-space
$H=\sqrt{2}h$).  Now, we observe that the first term in equation
\reff{discpde} is a discretization of minus the Laplacian of $r$ (in rotated
form) and that the second expression corresponds to minus the second
derivative of $r$ along the diagonal $x=y$ (which can be written as
$-1/2\Delta r -\partial^2 r/\partial x\partial y$). The right-hand-side of
\reff{discpde} approximates the function $3\sqrt{2}(6x(1-x)-1)\delta (x-y)$
over the diagonal. All together we obtain equation \reff{pde}, which shall be
interpreted in the weak sense.

To compute the solution of equation \reff{tij} we have developed a multigrid
method composed by red-black Gauss-Seidel relaxations, full-weighting of
residuals for fine-to-coarse mesh transfers and of bilinear interpolation for
coarse-to-fine corrections (see eg.  Hackbusch (1985) for definitions).  The
scheme is very fast, it takes only a few seconds to solve equation \reff{tij}
for $N=511$ on a Sun-workstation and around 15 minutes for the solution for
$N=4095$ on a Dec-alpha (the equation in this case has more than 16 million
unknowns). The equations were solved to the level of machine precision, using
double precision (16-byte) words.

\begin{figure}[ht]
\begin{center}
\epsfig{figure=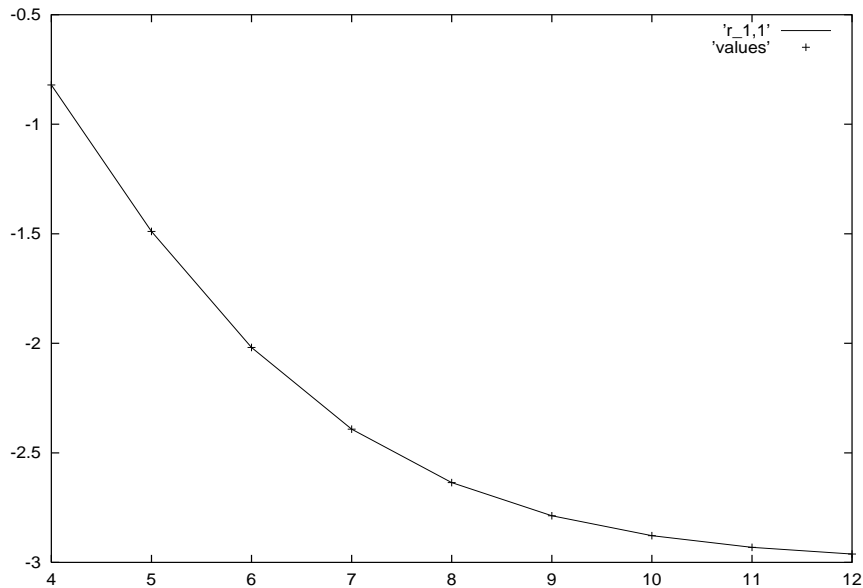,height=12cm,width=8cm,angle=-90}
\caption{ Graphic of the values of $R(1,1)/(N+1)^2$ 
for several values of $\log_2 (N+1)$.}
\end{center}
\end{figure}
\par\noindent
The results shown in Figures 1,2,3 were obtained with $\alpha=0$, but the same
asymptotics are expected for other values of $\alpha$.
In Figure 1 we present the function $r=R/(N+1)^3$, with a resolution
of $N=127$ (the function doesn't change much visually for larger
values of $N$).  In Figure 2 we plot $r=R/(N+1)^3$ over the diagonal
$\{i=j\}$, for several values of $N$ up to 511 - for larger values the graphics 
are indistinguishable. We can observe the convergence of
$r$ as $N$ grows (to the weak solution of the partial differential equation \reff{pde}).
The upper and lower constant bounds for $R/(N+1)^3$ we proved to exist in
theorem \reff{z6}, are attained over the diagonal $\{i=j\}$ and their
approximate values can be read from Figure 2.
\par\noindent
Another aspect we tried to access from the numerical solutions is the
behaviour of $R(1,1)/(N+1)^2$ (which corresponds to the derivative of $r$ over
the diagonal - the function of Figure 2 - at $x=0$). We computed these values
for $N=15$ up to $N=4095$ (with $N+1$ as powers of two).  The results are
presented graphically in Figure 3, where $N$ is displayed in $\log$ scale.
The results suggest that as $N\to\infty$ $R(1,1)/(N+1)^2$ converges to a
constant (around $-3$). This together with the definition of $\sigma^N$ given
in \reff{r} supports \reff{p2}. The constant $c(1,1)$ seems to be around $4$.

\section{The infinite silo model}

To prove Theorem \ref{111} we will use a spatially infinite limit of the
silo
model. This model represents the transmission of weight among weightless
grains in an infinitely large silo. For each $i\in\Z$ let $\eta_t(i)$
denote
the weight supported by site $i$ at time $t\ge 0$; the dynamics is
represented
by the system of equations: for $i \in \Z$ and $t\ge 1$ such that $it$ is
even,
\begin{equation}
\label{rap}
\eta_t(i) = \eta_{t-1}(i-1)U_{t-1}(i-1) + \eta_{t-1}(i+1) [
1-U_{t-1}(i+1)]\,,
\end{equation}
where $(U_t(i): i\in \Z, t\ge 1)$ are i.i.d.  uniform random variables in
$[0,1]$. We call this process \emph{infinite silo model} (ISM). Let
$\nu_{\rho}$ be a product measure on $[0,\infty)^{\Z}$ with marginals
Gamma$(2, 1/(2 \rho))$; $\rho$ is the expected value of each
marginal. Coppersmith et al. proved that the $\nu_{\rho}$ is invariant for
the process with periodic boundary conditions. We show now that $\nu_{\rho}$
is reversible (and hence invariant) for the infinite system. 

\begin{propos}
\label{pa2}  For each $\rho\ge 0$ the measure $\nu_\rho$ is reversible
for the ISM, that is, for all cylinder functions $f$ and $g$, $\E_{\nu_\rho}(
f(\eta_0) g(\eta_t)) = \E_{\nu_\rho}( g(\eta_0) f(\eta_t))$.
\end{propos}

The proof is based on the following
elementary lemma.
\begin{lemma} \label{lemma1} 
  Let $W$ be a random variable with distribution Gamma $(2,\rho/2)$.  Let $U$
  be a uniform random variable in $[0,1]$ independent of $W$. Then, $UW$ and
  $(1-U)W$ are independent exponentially distributed random variables with
  mean $\rho/2$.
\end{lemma}

\noindent{\bf Proof of Proposition \ref{pa2}}
It suffices to prove that if $\eta_0$ has law $\nu_\rho$  
then, for cylinder functions $f$ and $g$, 
\[ 
\E f(\eta_0)g(\eta_1) = \E g(\eta_0) f(\eta_1)\,. 
\]
Let
\[
X(i) = U_0(i) \eta_0(i)\,; \qquad  Y(i) = (1-U_0(i)) \eta_0(i)\,.
\]
Under $\nu_\rho$, $(\eta_0(i):i\in\Z)$ is a family of iid gammas random
variables. Then, by the lemma, $(X(i):i\in \Z)$ and $(Y(i):i\in \Z)$ are
independent families of iid exponential random variables of mean $\rho/2$ and
\begin{equation}
  \label{pa3}
 \E f(\eta_0)g(\eta_1) = \E f(X + Y) g(X'+Y')\,, 
\end{equation}
where $(X'(i),Y'(i)):= H(X,Y)(i) := (Y(i-1), X(i+1))$. Since $H(H(X,Y))=(X,Y)$
and $(X',Y')$ has the same law as $(X,Y)$ (iid exponential random variables),
\reff{pa3} equals to
\[
\E f(X'+Y') g(X + Y) = \E g(\eta_0) f(\eta_1)\,.\qquad\square
\]

\vskip 3mm

One of the tools to prove Theorem \ref{111} is to show that all invariant
and
translation invariant measures for the ISM can be written as convex
combinations of those product measures.

For any two initial configurations $\eta$ and $\xi$, we couple the two
processes $\eta_t$ and $\xi_t$ using the same uniform random variables to
update the process. Notice that this coupling is attractive, that is, 
\begin{equation}
  \label{z40}
  \eta(i) \le \xi(i) \mbox{ for all }i \in \Z, \mbox{ then }\eta_t(i) \le
  \xi_t(i), \mbox{ for all }i \in \Z\,. 
\end{equation}
Through out this section $\I$ and $\T$ ($\II$ and $\TT$) will
denote, respectively, the set of invariant and translation invariant
measures
for the process (coupled process) defined by \reff{rap}.

\begin{propos} 
\label{lemma2} 
If $\nu \in \II \cap \TT$, then 
\begin{equation}
\label{112}
\nu\{(\eta,\xi); \eta \ge \xi \mbox{ or } \eta \le \xi\} = 1\,.
\end{equation}
\end{propos}

The proof of this proposition is based on the following lemma.

\begin{lemma}
  \label{z50}
  Let $X$, $Y$ be identically distributed real random variables defined in the
  same probability space. Let $U$ and $V$ be identically distributed and
  independent random variables in $[0,1]$ with $\E U=1/2$; $U$ and $V$ are
  assumed to be independent of $X$ and $Y$.  If $Z= UX+VY$ is distributed as
  $X$ (and $Y$), then $\P(XY\ge 0)=1$.
\end{lemma}
\proof 
{From} the definition of $Z$,
\begin{eqnarray}
  \label{z51}
  \E(Z\one\{Z\ge 0\}) &=& \E(UX\one\{Z\ge 0\})+\E(VY\one\{Z\ge
  0\})\nonumber \\
&=& 2\E(UX\one\{Z\ge 0\})\nonumber \\
&=& 2\E(UX\one\{Z\ge 0, X\ge 0 \}) + 2\E(UX\one\{Z\ge 0, X< 0
  \})\nonumber \\
&=& 2\E(UX\one\{X\ge 0\})- 2\E(UX\one\{X\ge 0,Z< 0 \})  \})\nonumber \\
&&\qquad\qquad\qquad+\; 2\E(UX\one\{Z\ge 0, X< 0
  \})\,.\nonumber 
\end{eqnarray}
The second identity follows from the fact that $UX\one\{Z\ge 0\}$ and
$VY\one\{Z\ge 0\}$ have the same law. The third and fourth identities come
from set operations. Since $X$ and $Z$ are identically distributed and $U$ is
independent of $X$, $\E(Z\one\{Z\ge 0\}) = 2\E(UX\one\{X\ge 0\})$.
Canceling those terms we get
\[
0\;=\; -\,\E(UX\one\{X\ge 0,Z< 0 \})  \})\, + \,\E(UX\one\{Z\ge 0, X< 0
  \})\,.
\]
This implies 
\begin{eqnarray}
  \label{z60}
  \P(X> 0,Z< 0) &=& \P(X< 0,Z\ge 0) \;=\; 0\nonumber\\
\P(Y> 0,Z< 0) &=& \P(Y< 0,Z\ge 0) \;=\; 0
\end{eqnarray}
because $(X,Z)$ and $(Y,Z)$ have the same distribution. Hence
\[
\P(X> 0,Y< 0)\; =\; \P(X> 0,Y< 0, Z\ge 0)\,+\, \P(X> 0,Y< 0, Z< 0)\;=\; 0\,.
\]
The same argument shows that $\P(X< 0,Y> 0) =0$. 
\square

\paragraph{Proof of Proposition \ref{lemma2}}
Let $(\eta_0,\xi_0)$ be distributed according to $\nu \in \II\cap \TT$.
Then
the variables 
\begin{equation}
  \label{z55}
  X= \xi_0(i-1) - \eta_0(i-1);\quad Y = \xi_0(i+1) - \eta_0(i+1);\quad Z =
  \xi_1(i) - \eta_1(i) 
\end{equation}
satisfy the hypothesis of Lemma \ref{z50}. Hence, for all $i\in\Z$, 
\[
  \nu(\xi(i-1)=\eta(i-1)) = 1 \mbox{ or }
  \nu([\xi_0(i-1)-\eta_0(i-1)]\,[\xi_0(i+1)-\eta_0(i+1)]>0)=1\,,\nonumber  
\]
which implies the proposition.
\square

\begin{propos} \label{lemma3}
If $\mu \in (\I \cap \T)$ then
\begin{equation}
\label{119}
\mu = \int \nu_{\rho} \, d\lambda(\rho)\,,
\end{equation}
where $\lambda$ is a probability measure on $(0,1)$. 
\end{propos}
\proof For each $\rho\ge 0$ we can construct $\bar{\nu}_{\rho}\in\bar \I\cap
\bar \T$, a measure on $[0,\infty)^\Z \times [0,\infty)^\Z$ with marginals
$\mu$ and $\nu_{\rho}$ (Liggett (1976)).

We claim that $\mu \in \I\cap\T$ concentrate its mass on
configurations with asymptotic density. That is, for a configuration
$\eta$ define
\begin{equation}
\label{120}
M_n (\eta) = \frac{1}{2n +1} \sum_{i= -n}^{n} \eta(i)\,.
\end{equation}
Then, for all $a\in[0,1]$,
\begin{equation}
\label{121}
\mu \{\eta; \liminf_{n\to\infty} M_n(\eta) < a < \limsup_{n\to\infty}
 M_n(\eta)\} = 0\,.
\end{equation}
In fact, if \reff{121} is false, there exists an $a$ such that
\begin{equation}
\label{122}
\bar{\nu}_{a} \{(\eta, \xi); \eta \ge \xi \mbox{ or } \eta \le \xi\}^c
\ge \mu \{\eta; \lim \inf M_n(\eta) < a < \lim \sup M_n(\eta)\} > 0\,,
\end{equation}
which contradicts Proposition \ref{lemma2}. Therefore, taking $0 =
\rho_0 < \rho_1 < \ldots < \rho_m = 1$, with $\rho_i=i/m$, with $m$
positive integer, we have
\begin{equation}
\label{123}
\mu(\cdot) = \sum_{i=1}^{m} \mu(\cdot|C_i) \, p_i\,,
\end{equation}
where $C_i := \{ \eta; \lim_{n} M_n(\eta) \in (\rho_{i-1},\rho_{i}]\}$ and
$p_i = \mu(C_i)$, for $i = 1, \ldots, m$.  

Define
\begin{equation}
\label{123a}
\underline\mu_m(\cdot) = \sum_{i=1}^{m} \nu_{\rho_{i-1}} \, p_i\,;\qquad\qquad
\overline\mu_m(\cdot) = \sum_{i=1}^{m} \nu_{\rho_i} \, p_i\,.
\end{equation}
By definitions \reff{123} and \reff{123a} and Proposition \reff{lemma2} we can
couple $\mu$, $\underline\mu$ and $\overline\mu$ and show
$$
\underline\mu_m\; \le_{st}\; \mu \;\le_{st}\; \overline\mu_m\,.
$$
for all $m$, where $\le_{st}$ is the usual stochastic domination of
measures defined by: $\mu \le_{st} \nu$ if and only if for all non decreasing
function $f$, $\mu f \le \mu' f$; equivalently, $\mu \le_{st} \mu'$ if and
only if there exists a joint realization $(W,W')$ with marginals $\mu$ and
$\mu'$ such that $W\le W'$. Defining $\lambda(\rho) = \mu(\eta: \lim_{n}
M_n(\eta)\le \rho)$ we have
$$
\lim_{m\to\infty} \underline\mu_m = \lim_{m\to\infty} \overline\mu_m =
\int \nu_\rho\, d\lambda(\rho)
$$
This proves the proposition.  \hfill \square

\section{Proof of Theorem \ref{111}}

In this section we use the characterization of the set of invariant and
translation invariant measures for the ISM given by Proposition \ref{lemma3}
and the bounds on the covariances \reff{z7} to prove Theorem \ref{111}.  Let
$\mu^N$ be the unique invariant measure for the silo in the finite box $\lnt$.
We extend the definition of the measure $\mu^N$ in $[0,\infty)^\Z$ by setting
$\mu^N(\eta:\eta(i)=0) = 1$ for all $i \not\in \lnt$.

Let $U_t(i)$ be the uniform random variables defined in the introduction and
define $\mathbf U_t$ as the matrix with entries
\begin{eqnarray}
  \label{c8}
 \mathbf U_t(i,j) &:=& \cases {U_t(i) &if $j=i-1$\cr
1-U_t(i) &if $j=i+1$\cr
0&if $|j-i|>1$\,. }
\end{eqnarray}
Define the truncations
\begin{eqnarray}
 \label{cc8}
 [\mathbf U]^N(i,j) &:=& \mathbf U(i,j)\one\{i\in\Lambda^N\}\\
{[V]}^N(i) &:=& V(i)\one\{i\in\Lambda^N\}\,.\nonumber
\end{eqnarray}
The evolution equation \reff{W(i,t)} is then a product of a left
vector and a matrix plus an independent vector:
\begin{eqnarray}
  \label{c2}
  W^N_{t+1} &=& W^N_{t}[\mathbf U_t]^N + [V_{t+1}]^N\,.\nonumber
\end{eqnarray}
Notice that the definition of
$[\mathbf U_t]^N$ guarantees $W^N_t(j) \equiv 0$ for $j\notin\Lambda^N$.

\paragraph{Existence of the limit}
Let $W^N$ be a configuration distributed according to the invariant measure
$\mu^N$.  
The maximum of \reff{rw(i)} is attained for $i\in\{N/2,(N+1)/2\}$, hence 
\begin{equation}
\label{ca:100}
\E \Bigl[ \frac{W^N([rN]+ \ell)}{N^2} \Bigr] \;\le\; \frac14
\,+\,O\Bigl(\frac1{2N}\Bigr) \,,
\end{equation}
for all $\ell\in\Z$. Therefore, the sequence of probability
measures $\{\tau_{[rN]}\Theta_{N^2}\mu^N, N \ge 1\}$ (that is, $\{$law of
$\frac{W^N([rN]+\cdot)}{N^2}, N \ge 1\}$) is tight. Consider a convergent
subsequence $\{N_k, k \ge 1\}$ and let $\mu$ be the weak limit of this
subsequence:
\begin{eqnarray}
\label{c6}
\mu := \lim_{k\to\infty} \tau_{[rN_k]}\Theta_{N_k^2}\mu^{N_k}\,. 
\end{eqnarray}

\paragraph{Invariance of the limit}
The equation $\mu^N = S^N(t)\mu^N$ is equivalent to
\begin{equation}
  \label{c1}
  \E f(W^N) = \E f(W^N[\mathbf U]^N + [V]^N)\,,
\end{equation}
for cylinder bounded $f$, where $\mathbf U$ has the same law as $\mathbf U_t$,
$V$ is a vector with the same law as $V_{t+1}$ and $W^N$, $\mathbf
U$ and $V$ are independent.

If $\eta$ is distributed according to $\nu$, then the invariance of $\nu$ for
the infinite silo model is equivalent to 
\begin{equation}
  \label{c4}
  \E f(\eta) = \E f(\eta \mathbf U)\,,
\end{equation}
for cylinder bounded $f$, where $\eta$ and $\mathbf U$ are independent.

Let $\tilde\eta^k = \tau_{[rN_k]}\Theta_{N_k^2} W^{N_k}$, $\tilde V^k =
\tau_{[rN_k]} V$ and $\widetilde{\mathbf U}^k=\tau_{[rN_k]}\mathbf U$, where
$(\tau_\ell \mathbf U)(i,j) = \mathbf U(i-\ell,j-\ell)$.  Take $k$ so large such
that the support of $\tau_{[rN_k]}f$ is inside $\Lambda^{N_k}\setminus
\{1,N_k\}$ (and the truncations are unnecessary). Then the law of
$\tilde\eta^k$ must satisfy
\begin{equation}
  \label{c5}
  \E f(\tilde\eta^k) = \E f(\tilde\eta^k \widetilde{\mathbf U}^k +
  (\tilde V^k/N_k^2)) \,.
\end{equation}
By translation invariance, $\widetilde{\mathbf U}^k$ and $\mathbf U$ have the same
law, $\tilde V^k$ and $V$ have the same law and by \reff{c6} $\tilde\eta^k$
converges in distribution to $\mu$ as $k\to\infty$; hence $\mu$ must satisfy
\reff{c4}, that is, $\mu \in \I$.

\paragraph{Translation invariance of the limit}
Let $(W^{N-2}_0,W^N_0, W^{N+2}_0)$ be identically null and 
\begin{eqnarray}
  \label{c88}
\lefteqn{ (W^{N-2}_{t+1},W^N_{t+1}, W^{N+2}_{t+1})} \\
&=&
(W^{N-2}_{t}[\tau_1\mathbf U_{t}]^{N-2} + [\tau_1V_{t+1}]^{N-2},
W^{N}_{t}[\tau_1\mathbf U_{t}]^N + [\tau_1V_{t+1}]^N, 
W^{N+2}_{t}[\mathbf U_{t}]^{N+2}+ [V_{t+1}]^{N+2} )  \nonumber
\end{eqnarray}
for $t\ge 0$. Then $W^{N-2}_t\le W^N_t\le W^{N+2}_t$ for all $t$. By Theorem
\ref{2.4} $(W^{N-2}_t,W^N_t, W^{N+2}_t)$ converges in distribution to
$(W^{N-2},W^N, W^{N+2})$ a vector with marginals $\tau_1 \mu^{N-2}$, $\tau_1
\mu^N$ and $\mu^{N+2}$
also satisfying $W^{N-2}\le_{st} W^N\le_{st} W^{N+2}$, where $X\le_{st}Y$ means
that the law of $X$ is stochastically dominated by the law of $Y$ (see the
definition at the end of Section 4). This implies
\begin{equation} 
\label{130}
\tau_1 \mu^{N-2} \le_{st} \tau_1 \mu^N \le_{st} \mu^{N+2}\,.
\end{equation}
Analogously, the coupling 
\begin{eqnarray}
  \label{c888}
\lefteqn{ (W^{N-2}_{t+1},W^N_{t+1}, W^{N+2}_{t+1})} \\
&=&
(W^{N-2}_{t}[\tau_1\mathbf U_{t}]^{N-2} + [\tau_1V_{t+1}]^{N-2},
W^{N}_{t}[\mathbf U_{t}]^N + [V_{t+1}]^N, 
W^{N+2}_{t}[\mathbf U_{t}]^{N+2}+ [V_{t+1}]^{N+2} )  \nonumber
\end{eqnarray}
shows that
\begin{equation}
\label{130a}
\tau_1 \mu^{N-2} \le_{st} \mu^N \le_{st} \mu^{N+2}\,.
\end{equation}
{From} the coupling \reff{c888} and by Chebyshev,
\begin{eqnarray}
\P \Bigl( \frac{W^{N+2}(i) - W^N(i)}{N^2} \, > \, \epsilon
\Bigr) &\le& \frac{i(N+3-i) - i(N+1-i)}{N^2 \, \epsilon^2} 
\;=\; \frac{2i}{N^2 \, \epsilon^2} \rightarrow 0\nonumber\\
\P \Bigl( \frac{W^N(i)-W^{N-2}(i-1)}{N^2} \, > \, \epsilon
\Bigr) &\le& \frac{i(N+1-i) - (i-1)(N-i)}{N^2 \, \epsilon^2} 
\;=\; \frac{N}{N^2 \, \epsilon^2} \rightarrow 0\nonumber
\end{eqnarray}
as $N \rightarrow \infty$ for all $i \in \{0, 1, \ldots\}$.  Therefore, if
$\mu^{N_k} \tau_{rN_k}\Theta_{{N_k}^2}$ converges weakly to $\mu$ the same
is
true for $ \mu^{N_k-2} \tau_{rN_k+1}\Theta_{{N_k}^2}$ and $\mu^{N_k+2}
\tau_{rN_k}\Theta_{{N_k}^2}$. This and \reff{130} imply $\tau_1\mu=\mu$
and
hence that $\mu\in\cal T$.

\paragraph{Correlations of the limit.}
Since $\mu\in\I\cap\T$, by Proposition \ref{lemma3} $\mu$ is a convex
combination of gamma distributions with mean $r(1-r)$. The covariances of
$\mu$ equal the variance of the density of $\mu$: for all $i\neq j$,
\begin{eqnarray}
  \label{z4}
  \int \mu(d\eta) \eta(i)\eta(j) - \Bigl(\int \mu(d\eta) \eta(i)\Bigr)^2 
&=& \int \int  \nu_\rho(d\eta) \eta(i)\eta(j)\lambda(d\rho)
-\Bigl(\int \rho \lambda(d\rho)\Bigr)^2 \nonumber\\ 
&=& \int \rho^2 \lambda(d\rho)-\Bigl(\int \rho
\lambda(d\rho)\Bigr)^2\;\ge\;0\,.
\end{eqnarray}

Since the joint distribution $Q_k$ of the couple $({\eta(r{N_k}+i)\over
  {N_k}^2}, {\eta(r{N_k}+j)\over{N_k}^2})$ (under $\mu^{N_k}$) converges in
distribution to the joint distribution $Q$ of $(\eta(i),\eta(j))$ (under
$\mu$), then, as in the discussion following \reff{tig}, by Fatou's Lemma and
Skorohod representation for distribution convergence, $\int xy\, dQ \le
\liminf_{k\to\infty} \int xy\, Q_k(x,y)$.  Since by \reff{ca:100} the means
converge and $\lim_{k\to\infty}{\sigma^{N_k}(i,j) \over{N_k}^4}=0$ (by Theorem
\ref{z6}) the covariances of $\mu$ are nonpositive.

Conclusion: We have proved that $\mu$ is a convex combination of $\nu_\rho$
with mean $r(1-r)$ and with zero correlations. This means that $\mu=
\nu_{r(1-r)}$. Since each convergent subsequence of
$\{\mu^N\tau_{rN}\Theta_{N^2}, N \ge 1\}$ converges to $\nu_{r(1-r)}$, so does
the sequence. \square

\bigskip

\noindent{\bf Acknowledgments.} 
We thank Raul Rechtman for calling our attention to the silo with absorbing
walls, to Joachim Krug for insisting that the correlations were computable and
for mentioning relevant references and to Fabio Machado and Stefano Olla for
fruitful discussions. We also thank the referees for their comments and
suggestions. This work started when PAF and NLG visited Departamento de
Ingenieria Matem\'atica de la Universidad de Chile. We thank support from
FONDAP program in Stochastic Modeling, C\'atedra Presidencial, Funda\c c\~oes
Vitae-Antorchas-Andes, CNPq, PRONEX and FAPESP.

\parindent 0pt
Saulo R. M. Barros, Pablo A. Ferrari \\
IME USP, Caixa Postal 66281, 05315-970 - S\~{a}o Paulo, SP - BRAZIL - \\
email: {\tt saulo@ime.usp.br}, {\tt pablo@ime.usp.br}

\vskip 3mm
Nancy L. Garcia\\
IMECC, UNICAMP - Caixa Postal 6065, 13081-970 - Campinas, SP - BRAZIL -\\
email: {\tt nancy@ime.unicamp.br}
\vskip 3mm

Servet Mart\'{\i}nez\\
Departamento de Ingenier\'{\i}a Matem\'atica -
Casilla 170 - 3, Correo 3, Santiago, CHILE - \\
email: {\tt smartine@dim.uchile.cl}


\begin{thebibliography}{999}

  
\bibitem{claudin} Claudin P., Bouchaud J.P., Cates M.E., Wittmer J.P. (2000)
  Models of stress fluctuations in granular media {\sl Phys. Rev. E \bf
    57}: (4) 4441-4457.
  
\bibitem{coppersmith}
Coppersmith, S. N., Liu, C.-h., Majumdar, S., Narayan, O., Witten,
T. A. (1996) Model for force fluctuations in bead packs. {\sl Phys. Rev. E \bf
  53} 5:4673--4685. 

\bibitem{hack}
Hackbusch, W. (1985) {\sl Multi-grid methods and applications}, Springer.

\bibitem{harr}
Harr, M.E. (1977) {\sl Mechanics of Particulate Media}, McGraw-Hill.

\bibitem{kg} Krug J., Garcia J. (2000) 
     Asymmetric particle systems on R.
     {\sl J. Stat. Phys. \bf 99} (1-2) 31-55.


\bibitem{Le} Lewandowska, M., Mathur,   H., Yu, Y.-K. (2000) 
Dynamics and Critical Behaviour of the q-model. Preprint cond-mat/0007109.



\bibitem{L} Liggett, T. M. (1976) Coupling the simple exclusion process. {\sl
    Ann. Probab. \bf 4} 3:339--356.
  
\bibitem{liu} Liu, C. H., Nagel S. R., Schecter D. A., Coppersmith S. N.,
  Majumdar, S., Narayan O., Witten T. A. (1995) Force fluctuations in bead
  packs. {\sl Science \bf 269} 513--515.

\bibitem{peralta}
Peralta-Fabi R., M\'alaga C., Rechtman R., in Powders \& Grains, 
Behringer R. P., Jenkins J. T. Editors, A. A. Balkema 1997.

 
\bibitem{RM} Rajesh R., Majumdar S. N. (2000) Exact calculation of the
  spatiotemporal correlations in the Takayasu model and in the q model of
  force fluctuations in bead packs {\sl Phys. Rev. E \bf 62}: (3)
  3186-3196, Part A.
  
\bibitem{RMa} Rajesh R., Majumdar S. N. (2000a) 
     Conserved mass models and particle systems in one dimension
     {\sl J. Stat. Phys. \bf 99} (3-4) 943-965.


\bibitem{socolar} Socolar J. E. S. (1998) Average stresses and force fluctuations
  in noncohesive granular materials {\sl Phys. Rev. E \bf 57}: (3)
  3204-3215, Part B.

\end{thebibliography}
\end{document}